\documentclass{article}
\usepackage{here, amsmath, latexsym, amssymb, bm, ascmac, mathtools, multicol, tcolorbox, subfig, tikz, mathrsfs, siunitx, indentfirst}
\usepackage{amsthm}
\usepackage{tikz}
\usepackage{fancyhdr}
\usepackage[a4paper, total={6in, 9in}]{geometry}
\usepackage[T1]{fontenc}
\usepackage{lmodern}
\usepackage{booktabs,array,longtable}
\usepackage{enumitem,needspace}
\usepackage{xcolor}
\usetikzlibrary{arrows.meta,positioning,calc}
\usepackage[hidelinks]{hyperref}

\usepackage{titlesec}
\titleformat{\section}[block]{\bfseries\large}{\thesection.}{0.5em}{}
\titleformat{\subsection}[block]{\bfseries}{\thesubsection .}{0.5em}{\centering}

\theoremstyle{plain}
\newtheorem{thm}{Theorem}[section]
\newtheorem*{thm*}{Theorem}
\newtheorem{lem}[thm]{Lemma}

\newtheorem{prop}[thm]{Proposition}

\theoremstyle{definition}
\newtheorem{rem}{Remark}

\definecolor{navy}{HTML}{18344B}
\definecolor{blue}{HTML}{286791}
\definecolor{pale}{HTML}{EAF2F8}
\definecolor{teal}{HTML}{14776F}
\definecolor{paleteal}{HTML}{EAF5F1}
\definecolor{orange}{HTML}{B65B20}
\definecolor{paleorange}{HTML}{FCF0E5}
\definecolor{grey}{HTML}{56616B}
\definecolor{lightgrey}{HTML}{F2F4F5}

\DeclareMathOperator{\Tr}{Tr}
\DeclareMathOperator{\supp}{supp}
\DeclareMathOperator{\real}{Re}
\newcommand{\norm}[1]{\left\|#1\right\|}
\newcommand{\inpdt}[1]{\left<#1\right>}

\title{Critical orthonormal Strichartz estimates for Lebesgue spaces}
\author{Yinghao Gong}

\begin{document}
	\maketitle
	
	\begin{abstract}
		We study orthonormal Strichartz estimates at the critical coefficient exponent $q/2$ for fractional Schrödinger propagators. The proof for $d\geq3$ combines annular Schatten estimates with a frequency summation argument based on potential duality and Littlewood--Paley theory. In dimension two, a trace computation supplies the annular estimate needed for the global estimate. We also construct a counterexample to the global estimate with spatial $L^\infty$ norm and coefficients in $\ell^{q/2}$, even with a weak time norm.
	\end{abstract}
	
	\section{Introduction}
	
	Let $d\geq1$, $\mu\in\mathbb{R}\setminus\{0,1\}$. We denote the fractional Schrödinger propagators by \begin{equation*}
		e^{it|\nabla|^\mu}f(x):=\frac{1}{(2\pi)^d}\int_{\mathbb{R}^d}e^{i(x\cdot\xi+t|\xi|^\mu)}\widehat{f}(\xi)d\xi,
	\end{equation*}where $\widehat{f}$ is the Fourier transform of $f$. The \textit{orthonormal Strichartz estimates} for fractional Schrödinger propagators are the inequalities \begin{equation}\label{eq: ortho}
		\norm{\sum_j \lambda_j\left|e^{it|\nabla|^\mu}|\nabla|^{-s_\mu}f_j\right|^2}_{L_t^{q/2}L_x^{r/2}(\mathbb{R}\times\mathbb{R}^d)}\lesssim\norm{\lambda}_{\ell^{\alpha}},
	\end{equation}where $(f_j)_j\subset L^2(\mathbb{R}^d)$ is an orthonormal system, and $\lambda\in\ell^{\alpha}$. Here, \begin{equation}\label{eq: qr}
		2\leq q, r\leq\infty,\qquad\frac{2}{q}+\frac{d}{r}\leq\frac{d}{2},
	\end{equation}with $\alpha\geq1$, and\begin{equation}\label{eq: main scaling}
		s_\mu=\frac{d}{2}-\frac{d}{r}-\frac{\mu}{q}.
	\end{equation}Inequality \eqref{eq: ortho} fails if $s_\mu\neq d/2-d/r-\mu/q$, or \textcolor{blue}{if} $(q,r)$ does not satisfy \eqref{eq: qr}. We call the pair $(q,r)$ \textit{admissible} if it satisfies \eqref{eq: qr}. For a fixed $(q,r)$ and a fixed $\alpha_*$, we call inequality \eqref{eq: ortho} \textit{critical} with exponent $\alpha_*$ if \eqref{eq: ortho} holds when $\alpha=\alpha_*$, but fails if $\alpha>\alpha_*$. We refer to \eqref{eq: ortho} with $\alpha=\alpha_*$ as the critical estimate. For $\mu\in\mathbb{R}\setminus\{0,1\}$, \cite{Frank1,Frank3,Bez,Bez3} established critical orthonormal Strichartz estimates in the region below, with the additional assumption $d/r+\mu/q<d$ in \cite{Bez3}: \begin{equation*}
		d\geq2,\qquad2< q, r<\infty,\qquad\frac{1}{r}>\frac{d-1}{dq},\qquad\frac{2}{q}+\frac{d}{r}\leq\frac{d}{2}.
	\end{equation*}Moreover, for each admissible pair $(q,r)$ with $2<q,r<\infty$ with $1/r<(d-1)/dq$, \cite{Bez,Bez3} identified the necessary bound $\alpha\leq\alpha^*=q/2$ and left the critical case $\alpha=q/2$ unresolved in those works. 
	
	Orthonormal Strichartz estimates control the density formed by a family of solutions whose initial data are orthogonal in $L^2$. The coefficient norm records how much summability is gained from orthogonality beyond an estimate obtained by applying the scalar inequality to each function separately. In this paper we consider the orthonormal Strichartz estimates for fractional Schrödinger propagators and focus on the critical exponent. Our main purpose is to establish the estimate at this exponent in the following region. 
	\begin{thm}\label{thm: main}
		Let $d\geq3$, $\mu\in\mathbb R\setminus\{0,1\}$. Suppose that the pair $(q,r)$ satisfies \begin{equation}\label{geo: region}
			2< q, r<\infty,\qquad\frac{1}{r}<\frac{d-1}{dq},\qquad\frac{2}{q}+\frac{d}{r}<\frac{d}{2}.
		\end{equation}Then the inequality\begin{equation}\label{eq: main ortho}
			\norm{\sum_j \lambda_j\left|e^{it|\nabla|^\mu}|\nabla|^{-s_\mu}f_j\right|^2}_{L_t^{q/2}L_x^{r/2}}\lesssim\norm{\lambda}_{\ell^{q/2}}
		\end{equation}holds for every orthonormal system $(f_j)_j\subset L^2$ and every coefficient sequence $\lambda\in\ell^{q/2}$.
	\end{thm}
	
	For simplicity, we use the following notation and diagram. 
	
	\begin{figure}[H]
		\centering
		\begin{tikzpicture}
			\draw [->, dotted, thin] (0, 0) -- (6.3, 0) node [right] {$\frac{1}{r}$};
			\draw [->, dotted, thin] (0, 0) -- (0, 6.3) node [above] {$\frac{1}{q}$};
			\draw [loosely dotted] (6, 6) -- (2, 6) -- (6, 0) --cycle;
			\draw (0, 0) -- (0, 6) -- (2, 6) -- (6, 0) -- cycle;
			\draw [loosely dotted] (0, 3) -- (6, 6);
			\draw (0, 0) -- (3, 4.5);
			\draw [loosely dotted] (0, 3) -- (2, 3);
			\draw [loosely dotted] (0, 0) -- (2, 6);
			\fill (0, 0) circle (0.05cm) node [below left] {$O$};
			\fill (3, 4.5) circle (0.05cm) node [below] {$A$};
			\fill (6, 0) circle (0.05cm) node [below] {$B$};
			\fill (0, 6) circle (0.05cm) node [left] {$C$};
			\fill (2, 6) circle (0.05cm) node [above] {$D$};
			\fill (0, 3) circle (0.05cm) node [left] {$M$};
			\fill (2, 3) circle (0.05cm) node [right] {$N$};
		\end{tikzpicture}
		\caption{The exponent plane for $d=3$.}\label{fig:exponent-regions}
	\end{figure}
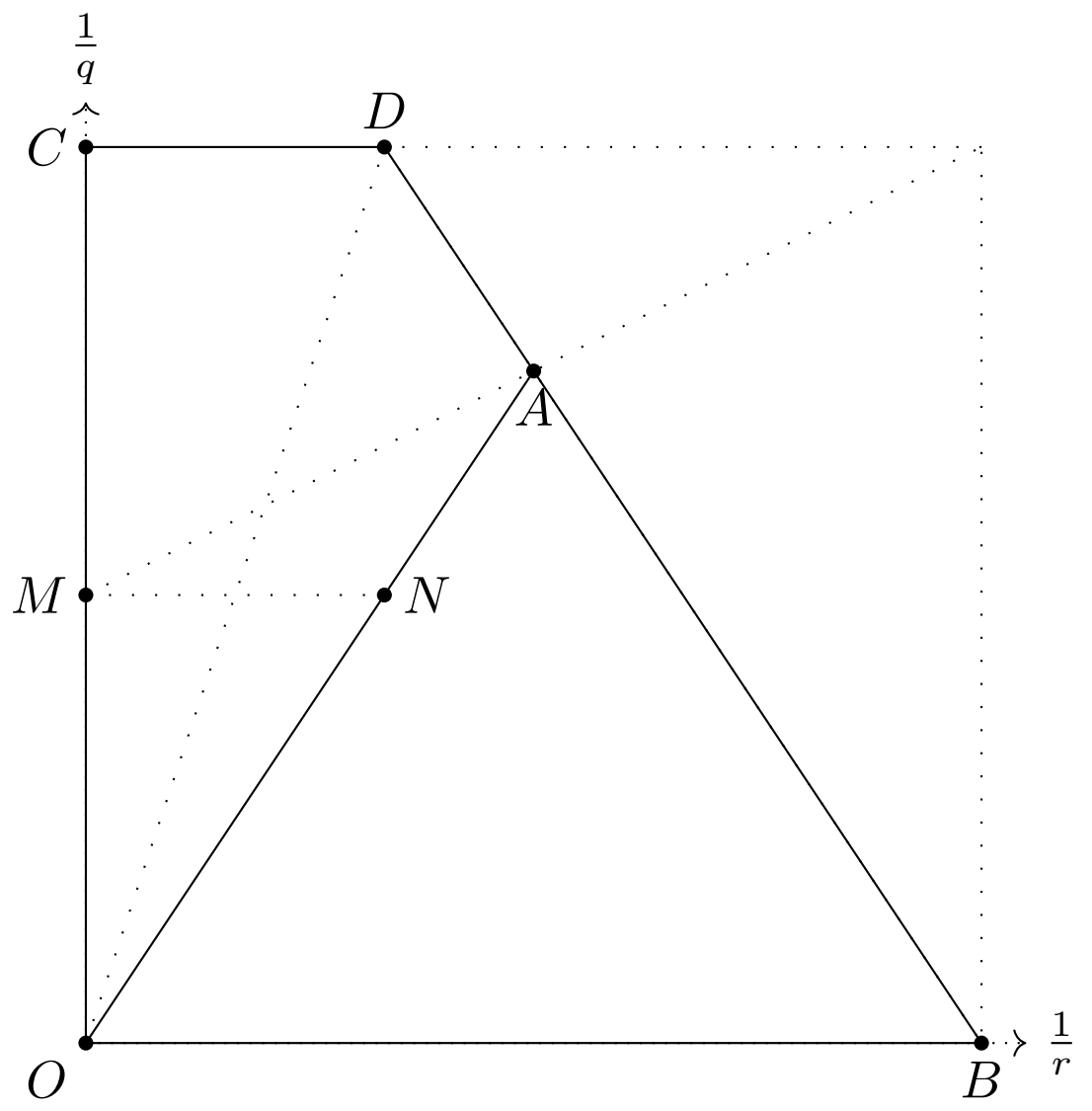
	
	In the \((1/r,1/q)\)-plane, set
	\begin{align*}
		O&=(0,0),
		&A&=\left(\frac{d-1}{2(d+1)},\frac{d}{2(d+1)}\right),
		&B&=\left(\frac12,0\right),
		&C&=\left(0,\frac12\right),\\
		D&=\left(\frac{d-2}{2d},\frac12\right),
		&M&=\left(0,\frac14\right),
		&N&=\left(\frac{d-1}{4d},\frac14\right).
	\end{align*}
	We use the usual notation for open, closed, and half-open line
	segments.
	For several points, \(X_1\cdots X_n\) denotes their convex hull,
	and \((X_1\cdots X_n)^o\) denotes its interior in \(\mathbb R^2\). The region of Theorem \ref{thm: main} and the critical segment are denoted respectively by $(OCDA)^o$ and $(O,A]$. 
	
	\vspace{2ex}
	We write $\mathcal{O}_\mu(q,r,\alpha)$ for the assertion that \eqref{eq: ortho} holds for every orthonormal system $(f_j)_j\subset L^2$ and every coefficient sequence $\lambda\in\ell^\alpha$, with $s_\mu$ satisfying \eqref{eq: main scaling}. By sequence embedding, $\mathcal{O}_\mu(q,r,\alpha_0)$ implies $\mathcal{O}_\mu(q,r,\alpha_1)$ if $1\leq\alpha_1<\alpha_0$. We also remark that the scalar classical Strichartz estimate
	\begin{equation}\label{eq: strichartz}
		\norm{e^{it|\nabla|^\mu}|\nabla|^{-s_\mu}f}_{L_t^qL_x^r}
		\lesssim \norm{f}_{L^2},
		\qquad \mu\notin\{0,1\},
	\end{equation}
	and the triangle inequality already give \eqref{eq: ortho} with $\alpha=1$. Thus, it is natural to ask whether the critical estimate $\mathcal{O}_\mu(q,r,\alpha_\ast)$ holds for each admissible pair.
	
	\vspace{2ex}
	\noindent\textbf{Frameworks:} Inequality \eqref{eq: strichartz} has been established by many authors as a useful tool to analyse the size and decay of solutions of dispersive partial differential equations. (See \cite{Strichartz,Ginibre-Velo,Keel-Tao}.) When $\mu=2$, this inequality was proved to hold for all admissible pairs, except at $(q,r)=(\infty,\infty)$ for $d\geq1$ and at $(q,r)=(2,\infty)$ for $d\geq2$. (See \cite{Montgomery-Smith,Z.Guo}.) For $\mu\notin\{0,1\}$ and finite admissible $(q,r)$, this inequality follows from annular dispersive estimates and Littlewood--Paley theory. (See \cite{Bez3}.) For $r=\infty$ and $\mu>0$, $\mu\ne1$, Guo--Li--Nakanishi--Yan \cite{Z.Guo} treat the admissible cases $2<q<\infty$ when $d\geq2$, and $4\leq q<\infty$ when $d=1$. Their Theorem~1.1 proves failure at $(q,r)=(2,\infty)$ for $d\geq3$ and $0<\mu\leq2$, $\mu\ne1$, but it does not settle this endpoint for $\mu>2$ or $\mu<0$. 
	
	For previous results on inequality \eqref{eq: ortho} with admissible pairs $(q,r)$ and $\mu=2$, we summarise the known results as follows: \begin{itemize}
		\item[1.] \textbf{On $[B,A]$:} Frank--Lewin--Lieb--Seiringer \cite{Frank1} first established $\mathcal{O}_2(q,r,\alpha_\ast)$ for every $d\geq1$, every $(q,r)$ satisfying $2<r\leq(2d+4)/d$ and $d/2=d/r+2/q$, with motivation from the study of systems of independent fermions. Their Proposition~2 also proves failure at $r=(2d+2)/(d-1)$ with $\alpha=(d+1)/d$, using a divergent sphere integral when $d\geq2$. Frank--Sabin \cite{Frank3} extended the positive range to $2<r<(2d+2)/(d-1)$ by using the Schatten duality principle.
		\item[2.] \textbf{On $(OAB)^o$ and $(OCDA)^o$:} For finite admissible exponents $2\leq q,r<\infty$, Bez--Hong--Lee-- Nakamura--Sawano \cite{Bez} established the optimal exponent $\alpha_\ast$ for $\mathcal{O}_2(q,r,\alpha)$ in the Sobolev range $d/2\geq d/r+2/q$: \begin{equation*}
			\alpha_*=\begin{cases}
				qrd/(2qd+2r),& 1/r>(d-1)/dq;\\ q/2,& 1/r\leq(d-1)/dq.
			\end{cases}
		\end{equation*}They proved $\mathcal{O}_2(q,r,\alpha_\ast)$ for every $1/r>(d-1)/dq$. They also proved that for $d\geq2$, $\mathcal{O}_2(q,r,\alpha)$ is true if $1/r\leq(d-1)/dq$, $d/2> d/r+2/q$, and $1\leq\alpha<\alpha_*=q/2$. For $q=\infty$, $2<r<\infty$, they claimed that $\mathcal{O}_2(q,r,\alpha_\ast)$ fails, by establishing a weak type failure of the corresponding velocity average estimates. 
		\item[3.] \textbf{On $[O,A]$:} \cite{Bez} also proved the failure of $\mathcal{O}_2(q,r,q/2)$ on the critical segment $1/r=(d-1)/dq$ for $d\geq2$. For $d=1$, the critical segment becomes $r=\infty$, $4\leq q\leq\infty$. Bez--Kinoshita--Shiraki \cite{Bez2} proved $\mathcal{O}_2(q,\infty,\alpha)$ for $4<q<\infty$, $1\leq\alpha<q/2$, and explained the failure of the critical restricted weak type estimate for $4\leq q<\infty$.
		\item[4.] \textbf{On $(O,C)$:} On $r=\infty$ with $d\geq2$, \cite{Bez2} also claimed $\mathcal{O}_2(q,r,\alpha)$ for $2<q<\infty$, $1\leq\alpha< q/2$. Moreover, they sketched restricted weak type estimates with $\ell^{q/2,1}$ coefficients and an $L_t^{q/2,\infty}L_x^\infty$ target.
		\item[5.] \textbf{On $[A,D]$:} On the sharp admissible segment $(2d+2)/(d-1)<r<2d/(d-2)$, $d\geq2$, \cite{Bez} provided weak type estimates for $\alpha=q/2$. Hadama--Toyoshima \cite{HT} proved Conjecture~1.3 in \cite{Bez} at $A$ for $d\geq2$, by a direct singular value counting method. For $d\geq3$, interpolating with the Keel--Tao endpoint estimate yields the strong estimate $\mathcal{O}_2(q,r,q/2)$ on the entire open segment. (See \cite[Corollary 1.6]{HT}.)
		\item[6.] \textbf{On $(M,N)$:} In the study of global well-posedness for an infinite nonlinear Schr\"odinger system, Chen--Hong--Pavlovi\'c \cite{Chen} proved density estimates with inhomogeneous Sobolev weights. The homogeneous adaptation and Sobolev embedding yield $\mathcal{O}_2(4,r,2)$ for $d\geq2$ and $4d/(d-1)<r<\infty$.
	\end{itemize}
	
	For fractional propagators with $\mu\notin\{0,1,2\}$, the known results are as follows: \begin{itemize}
		\item[1.] \textbf{On $(OAB)^o$ and $(OCDA)^o$:} Bez--Lee--Nakamura \cite{Bez3} extended the results of \cite{Bez} to fractional propagators, including negative $\mu$, under the additional condition $d/r+\mu/q<d$. For finite admissible exponents $2\leq q,r<\infty$ with $1/r>(d-1)/(dq)$, they proved $\mathcal{O}_\mu(q,r,\alpha_*)$. Their published non-sharp theorem treats $d\geq3$ and the strict interior with $1\leq\alpha<\alpha_*=q/2$. Feng--Lee--Lee \cite{FengLeeLee} gave the critical Schr\"odinger estimate throughout the interior in Theorem~\ref{thm: main}, whereas for their fractional result, they obtained $\alpha<q/2$.
		\item[2.] \textbf{On $(M,N)$:} For $d\geq2$, $q=4$ and $4d/(d-1)<r<\infty$, Gong \cite{Gong} proved $\mathcal{O}_\mu(q,r,q/2)$ for $\mu>0$, $\mu\ne1$, by Hilbert--Schmidt kernel estimates and bilinear real interpolation. \cite{Gong} also identified a specific real interpolation failure in the argument of \cite{Feng}.
	\end{itemize}
	
	\vspace{2ex}
	
	\noindent\textbf{Strategy:} The first step in the proof of Theorem \ref{thm: main} is providing an estimate on a fixed annulus. It gives a larger exponent region and allows an arbitrary regularity exponent. Let $Q$ be the Fourier multiplier given by $\widehat{Qf}(\xi)=\chi(\xi)\widehat{f}(\xi)$, where $\chi\in C_c^\infty$ is real-valued and supported away from the origin. 
	
	\begin{thm}\label{thm: annulus}
		Suppose $d\geq3$, $\mu\in\mathbb R\setminus\{0,1\}$, $s\in\mathbb{R}$. If $(1/r,1/q)\in OCDA\setminus \{A\}$, then the inequality\begin{equation}\label{eq: annulus OCDA}
			\norm{\sum_j\lambda_j\left|e^{it|\nabla|^{\mu}}|\nabla|^{-s}Qf_j\right|^2}_{L_t^{q/2}L_x^{r/2}}\lesssim_Q\norm{\lambda}_{\ell^{q/2}}
		\end{equation}holds for every orthonormal system $(f_j)_j\subset L^2$ and every coefficient sequence $\lambda\in\ell^{q/2}$.
	\end{thm}
	
	The proof of Theorem \ref{thm: annulus} on $OCDA\setminus[A,D)$ follows the annular Schatten argument of \cite[Section 3]{Bez}, using oscillatory integral theory. The abstract result of Hadama--Toyoshima \cite[Corollary 1.6]{HT} supplies $(A,D)$. We then pass to the global estimate by summing the frequency-localised potential operators and applying Littlewood--Paley theory directly on $(OCDA)^o$. We remark that this summation method avoids the interpolation difficulty in \cite{Bez}, which comes from the strong type obstruction on $[O,A]$.
	
	\vspace{2ex}
	
	On $(O,C)$, the global estimate with coefficients in $\ell^{q/2}$ fails even if its time norm is weakened. The counterexample uses an orthonormal system supported in disjoint frequency annuli, whose evolved density is large near the spatial origin. This conclusion is compatible with restricted weak estimates using $\ell^{q/2,1}$, as considered in \cite{Bez2}.
	
	When $d=2$, we consider the triangle $OAC$. We remark that the failure of the annulus estimate at $C$ for $\mu=2$ was proved in \cite{Montgomery-Smith,T.Tao 1} by constructing a counterexample from Brownian motion. We establish a potential duality argument for proving \eqref{eq: annulus OCDA} at points on $(O,C)$ with integer $q/(q-2)\geq2$, and then use interpolation to obtain every point on $(O,C)$. We then obtain the global estimate by adapting the proof of Theorem \ref{thm: main}.
	
	We also remark that on $(A,D)$, inequality \eqref{eq: main ortho} follows directly from \cite[Corollary 1.6]{HT} when $d\geq3$ and $\mu=2$. For $\mu\neq2$, $s_\mu=(2-\mu)/q$ on this segment, so $e^{it|\nabla|^\mu}|\nabla|^{-s_\mu}$ is not bounded on $L^2$, as required by that result. The global estimate \eqref{eq: main ortho} on $(A,D)$ for $d\geq2$ and $\mu\neq0,1,2$ is not established by the present argument. 
	
	\vspace{2ex}
	
	\textbf{Organisation.} Section 2 recalls some definitions and the duality principles. Section 3 proves the annular estimates, and Section 4 establishes the global estimate in higher dimensions. Section 5 gives the spatial endpoint counterexamples. Section 6 proves the two-dimensional estimates.
	
	\vspace{2ex}
	
	\textbf{Acknowledgements.} I am grateful to my supervisor Professor Neal Bez for his useful comments, patient guidance, and continuous support. I would like to thank my friends Sonae Hadama and Shunya Toyoshima for their valuable suggestions. 
	
	\section{Preliminaries}
	
	We recall that $A\lesssim B$ means $A\leq CB$, where $C$ is a positive constant. More generally, $A\lesssim_{a}B$ means $A\leq C_{a}B$ where $C_{a}$ is a positive constant depending only on $a$. 
	
	Define  \begin{equation*}
		\widehat{f}(\xi):=\int_{\mathbb{R}^d}e^{-ix\cdot\xi}f(x)dx,\qquad f(x)=\frac{1}{(2\pi)^d}\int_{\mathbb{R}^d}e^{ix\cdot\xi}\widehat{f}(\xi)d\xi.
	\end{equation*}
	
	\begin{prop}[Lemma 3, \cite{Frank3}; Proposition 2.1, \cite{Bez}; Remark 2.3, \cite{FMSW}]\label{prop: duality}
		Suppose $2<q<\infty$, $2\leq r\leq\infty$ and $p,\alpha\geq1$. Then, for a bounded operator \begin{equation*}
			A: L^2\longrightarrow L_t^{q,2p}L_x^{r},
		\end{equation*}the inequality \begin{equation*}
			\norm{\sum_j\lambda_j\left|Af_j\right|^2}_{L_t^{q/2,p}L_x^{r/2}}\lesssim\norm{\lambda}_{\ell^{q/2,\alpha}}
		\end{equation*}holds for every orthonormal system $(f_j)_j\subset L^2$ and every coefficient sequence $\lambda\in\ell^{q/2,\alpha}$, if and only if \begin{equation*}
			\norm{WAA^*\overline{W}}_{\mathfrak{S}^{(q/2)',\alpha'}}\lesssim\norm{W}_{L_t^{2(q/2)',2p'}L_x^{2(r/2)'}}^2
		\end{equation*}holds for all $W\in L_t^{2(q/2)',2p'}L_x^{2(r/2)'}$.
	\end{prop}
	
	We use the notation for \textit{Lorentz spaces} in \cite{Bergh,Bez}. For $1\leq p<\infty$, $\mathfrak{S}^p=\mathfrak{S}^p(L^2)$ denotes the \textit{Schatten space} on $L^2$, which is the set of all compact operators $A$ on $L^2$ such that $\norm{A}_{\mathfrak{S}^p}=\norm{(\sigma_j(A))_j}_{\ell^p}<\infty$, where $\sigma_j(A)$ are \textcolor{blue}{the} singular values of the operator $A$. Also, $\mathfrak{S}^\infty=\mathfrak{S}^\infty(L^2)$ denotes the bounded operator space $\mathcal{B}(L^2)$. We refer the reader to \cite{Bergh} for the definition of Lorentz spaces, and \cite{Simon} for the properties of Schatten spaces.
	
	\vspace{2ex}
	
	\noindent\textbf{Potential duality:} We introduce another duality principle which was used in \cite{Frank1,FMSW}. Let  $d\geq2$, $\mu\neq0,1$. Fix $s\in\mathbb{R}$. Write $U(t)=e^{it|\nabla|^\mu}|\nabla|^{-s}$. Let $(f_j)_j$ be  an orthonormal system in $L^2$ and define the following operator\begin{equation}\label{eq: gamma}
		\gamma:=\sum_{j\in\mathbb{Z}}\lambda_j |f_j\rangle\langle f_j|.
	\end{equation}For $k\in\mathbb{Z}$, we have $\gamma f_k=\lambda_k f_k$, so that $(\lambda_j)_{j\in\mathbb{Z}}$ are precisely the eigenvalues of $\gamma$ with respect to eigenvectors $(f_j)_{j\in\mathbb{Z}}$, which means $\norm{(\lambda_j)_{j\in\mathbb{Z}}}_{\ell^\alpha}=\norm{\gamma}_{\mathfrak{S}^\alpha}$. By using the definition of the density function of an operator, we can write \eqref{eq: annulus OCDA} as \begin{equation*}
		\norm{\rho_{\gamma_Q(t)}}_{L_t^{q/2}L_x^{r/2}}\lesssim\norm{\gamma}_{\mathfrak{S}^{q/2}},
	\end{equation*}where $\gamma_Q(t)=U(t)Q\gamma QU(t)^*$. On the other hand, since $Q$ is a Fourier multiplier supported away from the origin, using Fubini's theorem, trace cyclicity \cite[Theorem 3.1]{Simon} and \eqref{eq: gamma} yields \begin{align}
		\Tr \left(\gamma\int_{\mathbb{R}}U(t)^*QV(t,x)QU(t)dt\right)&=\sum_{j\in\mathbb{Z}}\lambda_j\inpdt{\left(\int_{\mathbb{R}}U(t)^*QV(t,\cdot)QU(t)dt\right)f_j,f_j}\notag\\&=\sum_{j\in\mathbb{Z}}\lambda_j\int_{\mathbb{R}}\inpdt{U(t)^*QV(t,\cdot)QU(t)f_j,f_j}dt\notag\\&=\int_{\mathbb{R}}\sum_{j\in\mathbb{Z}}\lambda_j\inpdt{U(t)^*QV(t,\cdot)QU(t)f_j,f_j}dt\notag\\&=\int_{\mathbb{R}}\Tr(V(t,\cdot)\gamma_Q(t))dt\notag\\&=\int_{\mathbb{R}^{d+1}}V(t,x)\rho_{\gamma_Q(t)}(x)dxdt.\label{eq: Vrho}
	\end{align}Thus, by applying the standard duality argument to Schatten spaces \cite[Theorem 3.2]{Simon}, we obtain that inequality \eqref{eq: annulus OCDA} is equivalent to \begin{equation}\label{eq: BvQ}
		\norm{B_{V,Q}^s}_{\mathfrak{S}^{(q/2)'}}\lesssim_Q \norm{V}_{L_t^{(q/2)'}L_x^{(r/2)'}},
	\end{equation}where \begin{equation*}
		B_{V,Q}^s:=\int_{\mathbb{R}}U(t)^*QV(t,x)QU(t)dt.
	\end{equation*}Furthermore, removing $Q$, we observe that inequality \eqref{eq: main ortho} is equivalent to \begin{equation}\label{eq: Bv}
		\norm{B_{V}^{s_\mu}}_{\mathfrak{S}^{(q/2)'}}\lesssim \norm{V}_{L_t^{(q/2)'}L_x^{(r/2)'}},
	\end{equation}where \begin{equation*}
		B_{V}^{s_\mu}:=\int_{\mathbb{R}}|\nabla|^{-s_\mu}e^{-it|\nabla|^\mu}V(t,x)e^{it|\nabla|^\mu}|\nabla|^{-s_\mu}dt.
	\end{equation*}
	
	\section{Annulus estimates}\label{sec: annulus}
	
	In this section, we prove \eqref{eq: annulus OCDA} on $OCDA\setminus[A,D)$, by an adaptation of \cite{Bez}. We also give a short proof of \eqref{eq: annulus OCDA} on $(A,D)$, as an application of \cite{HT}.
	
	\vspace{2ex}
	
	Fix $\mu\in\mathbb{R}\setminus\{0,1\}$. As in the proof of \cite[(3.1)]{Bez}, we use Bessel's inequality, Plancherel's theorem, and the assumption $0\notin\supp \chi$ to obtain \begin{align*}
		\sum_j|\lambda_j|\left|e^{it|\nabla|^\mu}|\nabla|^{-s}Qf_j\right|^2&\leq\left(\sup_j|\lambda_j|\right)\sum_j\left|e^{it|\nabla|^\mu}|\nabla|^{-s}Qf_j\right|^2\\&\lesssim\norm{\lambda}_{\ell^\infty}\norm{|\xi|^{-s}\chi}_{L_\xi^2}^2\\&\lesssim_\chi\norm{\lambda}_{\ell^\infty}.
	\end{align*}This proves \eqref{eq: annulus OCDA} at $O$.
	
	Let $d\geq3$. For $(1/r,1/q)\in [C,D]$, we first remark that the Hessian of $|\xi|^\mu$ is nondegenerate on $\supp\chi$ for $\mu\neq0,1$, so the argument in \cite[Chapter 8, Proposition 2.5]{M.Stein 1} gives\begin{equation}\label{eq: dispersive}
		\sup_x\left|\int_{\mathbb{R}^d}e^{i(x\cdot\xi+t|\xi|^\mu)}|\xi|^{-s}\chi(\xi)d\xi\right|\lesssim(1+|t|)^{-d/2}.
	\end{equation}This bound also holds with amplitude $|\xi|^{-2s}\chi(\xi)^2$, which we will consider later. On the other hand, Plancherel's theorem gives $\norm{e^{it|\nabla|^\mu}|\nabla|^{-s}Qf}_{L^2}\lesssim\norm{f}_{L^2}$. Thus, \cite[Theorem 1.2]{Keel-Tao} implies \begin{equation*}
		\norm{e^{it|\nabla|^\mu}|\nabla|^{-s}Qf}_{L_t^{q}L_x^{r}}\lesssim\norm{f}_{L^2}
	\end{equation*}for every $(1/r,1/q)\in [C,D]$. Therefore, the triangle inequality gives \begin{align*}
		\left\|
		\sum_j\lambda_j
		\bigl|e^{it|\nabla|^\mu}|\nabla|^{-s}Qf_j\bigr|^2
		\right\|_{L_t^1L_x^{r/2}}
		&\le
		\sum_j|\lambda_j|
		\left\|
		e^{it|\nabla|^\mu}|\nabla|^{-s}Qf_j
		\right\|_{L_t^2L_x^r}^{\,2}\\
		&\lesssim
		\sum_j|\lambda_j|\|f_j\|_{L^2}^2\\
		&=
		\|\lambda\|_{\ell^1}
	\end{align*}for every $r\geq 2d/(d-2)$. This is \eqref{eq: annulus OCDA} on $[C,D]$. 
	
	\vspace{2ex}
	
	\noindent\textbf{On the segment $(A,D)$:} Since $\chi$ is compactly supported away from the origin, $e^{it|\nabla|^\mu}|\nabla|^{-s}Q$ is bounded as an $L^2\to L^2$ operator. Thus, \eqref{eq: dispersive} shows that when $d\geq3$, $e^{it|\nabla|^\mu}|\nabla|^{-s}Q$ satisfies all assumptions in \cite[Corollary 1.6]{HT}, and therefore we obtain \eqref{eq: annulus OCDA} on $(A,D)$.
	
	\vspace{2ex}
	
	\noindent\textbf{In the region $OCDA\setminus[A,D)$:} Define \begin{equation*}
		TG(t)=\int_{\mathbb{R}}e^{i(t-t')|\nabla|^\mu}|\nabla|^{-2s}Q^2G(t')dt'.
	\end{equation*}By Proposition \ref{prop: duality} and complex interpolation, Theorem \ref{thm: annulus} on $OCDA\setminus[A,D)$ reduces to proving that\begin{equation}\label{eq: annulus NA}
		\norm{W_1TW_2}_{\mathfrak{S}^{(q/2)'}}\lesssim\norm{W_1}_{L_t^{2(q/2)'}L_x^{2(r/2)'}}\norm{W_2}_{L_t^{2(q/2)'}L_x^{2(r/2)'}}
	\end{equation}holds for every $W_1,W_2\in L_t^{2(q/2)'}L_x^{2(r/2)'}$, and $(1/r,1/q)\in[N,A)$. We follow \cite[Section 3]{Bez}, using \eqref{eq: dispersive} and the kernel bound below.
	
	Choose $\psi\in C_c^\infty(\mathbb{R})$ such that \begin{equation*}
		\supp\psi\subset\left\{t;\frac{1}{2}\leq|t|\leq2\right\},\qquad\widehat{\psi}(0)=0,\qquad\sum_{k\in\mathbb{Z}}\psi(2^{-k}t)=1\text{ for every }t\neq0.
	\end{equation*}Define \begin{equation*}
		T_kG(t)=\int_{\mathbb{R}}\psi(2^{-k}(t-t'))e^{i(t-t')|\nabla|^\mu}|\nabla|^{-2s}Q^2G(t')dt',
	\end{equation*}so that $T=\sum_kT_k$, and split $T_k=T_{k,0}+T_{k,1}$ according to \begin{equation*}
		|x-y|>C2^k,\qquad|x-y|\leq C2^k.
	\end{equation*}Choose $C= 10\sup_{\xi\in\supp\chi}|\nabla_\xi(|\xi|^\mu)|$ (say). Integration by parts gives, for every $N\geq0$, \begin{equation*}
		\left|\int_{\mathbb{R}^d}e^{i(x\cdot\xi+t|\xi|^\mu)}|\xi|^{-2s}\chi(\xi)^2d\xi\right|\lesssim(1+|x|)^{-N}, \qquad2|x|>C|t|.
	\end{equation*}
	
	Firstly, we show the Hilbert--Schmidt inequality \begin{equation}\label{eq: annulus HS}
		\norm{W_1\sum_{k\in\mathbb{Z}}2^{zk}T_kW_2}_{\mathfrak{S}^{2}}\lesssim\norm{W_1}_{L_t^{4}L_x^{v}}\norm{W_2}_{L_t^{4}L_x^{v}},
	\end{equation}where $2<v<4$, and $\real z=-\beta(v,v)$. Here, \begin{equation*}
		\beta(v_1,v_2):=\frac{d+1}{2}-d\left(\frac{1}{v_1}+\frac{1}{v_2}\right).
	\end{equation*}

	For $T_{k,0}$, choose a sufficiently large $N>0$ and apply Young's convolution inequality in space and time to obtain\begin{equation*}
		\norm{(W_1T_{k,0}W_2)_k}_{\ell_{-\beta(v,v)}^1(\mathfrak{S}^2)}\lesssim\sum_{k\in\mathbb{Z}}\frac{2^{(1/2-\beta(v,v))k}}{(1+2^k)^N}\norm{W_1}_{L_t^4L_x^v}\norm{W_2}_{L_t^4L_x^v}\lesssim\norm{W_1}_{L_t^4L_x^v}\norm{W_2}_{L_t^4L_x^v}
	\end{equation*}for every $2<v<4$, where $\ell_\beta^p(X)$ denotes the $\beta$-weighted $\ell^p$ space equipped with the norm $\norm{(T_k)_k}_{\ell_\beta^p(X)}^p:=\sum_k2^{kp\beta}\norm{T_k}_{X}^p$ for $p<\infty$, and $\norm{(T_k)_k}_{\ell_\beta^\infty(X)}:=\sup_k2^{k\beta}\norm{T_k}_X$. (See, \cite{Bergh}.) The triangle inequality then implies\begin{equation}\label{eq: Tk0}
		\norm{W_1\sum_{k\in\mathbb{Z}}2^{-\beta(v,v)k}T_{k,0}W_2}_{\mathfrak{S}^2}\lesssim\norm{W_1}_{L_t^4L_x^v}\norm{W_2}_{L_t^4L_x^v}.
	\end{equation}
	
	Now, fix $2<v_*<4$. For $T_{k,1}$, let $\mathcal{T}_{v_*}=(\mathcal{T}_{k,v_*})_k$ be a vector-valued bilinear operator, where\begin{align*}
		\mathcal{T}_{k,v_*}(V_1,V_2):=2^{-2\beta(v_*,v_*)k}\int_{\mathbb{R}^2}\int_{|x-y|\lesssim2^k}&|\psi(2^{-k}(t-t'))|V_1(t,x)V_2(t',y)\\&\left|\int_{\mathbb{R}^d}e^{i\left((x-y)\cdot\xi+(t-t')|\xi|^\mu\right)}|\xi|^{-2s}\chi(\xi)^2d\xi\right|^2dxdydtdt'.
	\end{align*}Note that the Cauchy--Schwarz inequality gives\begin{equation*}
		\norm{W_1\sum_{k\in\mathbb{Z}}2^{-\beta(v_*,v_*)k}T_{k,1}W_2}_{\mathfrak{S}^2}^2\lesssim\sum_{k\in\mathbb{Z}}\mathcal{T}_{k,v_*}\left(|W_1|^2,|W_2|^2\right).
	\end{equation*}Thus, we shall prove the boundedness of $\mathcal{T}_{v_*}:L_t^2L_x^{v_*/2}\times L_t^2L_x^{v_*/2}\longrightarrow\ell^1$ to obtain \begin{equation}\label{eq: Tk1}
		\norm{W_1\sum_{k\in\mathbb{Z}}2^{-\beta(v_*,v_*)k}T_{k,1}W_2}_{\mathfrak{S}^2}\lesssim\norm{W_1}_{L_t^4L_x^{v_*}}\norm{W_2}_{L_t^4L_x^{v_*}}.
	\end{equation}
	
	For $2<v_1,v_2<4$, Young's convolution inequality in space and time gives\begin{equation}\label{eq: Tkv}
		\left|\mathcal{T}_{k,v_*}(V_1,V_2)\right|\lesssim2^{2\left(\beta(v_1,v_2)-\beta(v_*,v_*)\right)k}\norm{V_1}_{L_t^2L_x^{v_1/2}}\norm{V_2}_{L_t^2L_x^{v_2/2}}.
	\end{equation}Note that the upper bound above is not summable when $v_1=v_2=v_*$. To gain room for summability, we choose nearby exponent pairs so that
	interpolation yields the target point $(v_1,v_2)=(v_*,v_*)$. Denote $\gamma(v_1,v_2):=2\beta(v_1,v_2)-2\beta(v_*,v_*)$. Choose $\delta>0$ small enough so that $2<v_0,v_1<4$, where\begin{equation*}
		\frac{1}{v_0}=\frac{1}{v_*}+\delta,\qquad\frac{1}{v_1}=\frac{1}{v_*}-2\delta.
	\end{equation*}This choice ensures $\gamma(v_0,v_0)<0<\gamma(v_0,v_1)=\gamma(v_1,v_0)$. Write\begin{equation*}
		A_0=L_t^2L_x^{v_0/2},\qquad A_1=L_t^2L_x^{v_1/2},\qquad B_0=\ell_{-\gamma(v_0,v_0)}^\infty,\qquad B_1=\ell_{-\gamma(v_0,v_1)}^\infty.
	\end{equation*}Inequality \eqref{eq: Tkv} gives the following boundedness properties:\begin{equation*}
		\mathcal{T}_{v_\ast}:\begin{cases}
			A_0\times A_0\longrightarrow B_0,\\
			A_0\times A_1\longrightarrow B_1,\\
			A_1\times A_0\longrightarrow B_1.
		\end{cases}
	\end{equation*}Therefore, by \cite[Section 3.13, Exercise 5b]{Bergh}, as applied in \cite[(3.15)]{Bez}, we obtain the boundedness of $\mathcal{T}_{v_*}:L_t^2L_x^{v_*/2}\times L_t^2L_x^{v_*/2}\longrightarrow\ell^1$. This proves \eqref{eq: Tk1}, and \eqref{eq: annulus HS} then follows by \eqref{eq: Tk0} and \eqref{eq: Tk1}. 
	
	Using $\widehat{\psi}(0)=0$, the argument in \cite[(3.17), (3.18)]{Bez} gives \begin{equation}\label{eq: annulus op}
		\norm{W_1\sum_{k\in\mathbb{Z}}2^{zk}T_kW_2}_{\mathfrak{S}^\infty}\lesssim\norm{W_1}_{L_t^\infty L_x^\infty}\norm{W_2}_{L_t^\infty L_x^\infty}
	\end{equation}for every $z$ satisfying $\real z=-1$. Therefore, using (analytic) complex interpolation yields \begin{equation*}
		\norm{W_1\sum_{k\in\mathbb{Z}}T_kW_2}_{\mathfrak{S}^{(q/2)'}}\lesssim\norm{W_1}_{L_t^{2(q/2)'} L_x^{2(r/2)'}}\norm{W_2}_{L_t^{2(q/2)'} L_x^{2(r/2)'}}
	\end{equation*}at $z=0$. Indeed, for every fixed $2<v_*<4d/(d+1)<4$, the above inequality follows from the
	interpolation identity \begin{equation*}
		\left(0,1-\frac{2}{q},\frac{1}{2}-\frac{1}{q},\frac{1}{2}-\frac{1}{r}\right)
		=(1-\theta)(-1,0,0,0)
		+\theta\left(
		-\beta(v_*,v_*),\frac12,\frac14,\frac1{v_*}
		\right),
	\end{equation*}where $0<\theta=1/(1-\beta(v_*,v_*))<1$. Moreover, we have \begin{equation*}
		\frac{1}{q}=\frac{d}{(d-1)r},
		\qquad
		\frac{1}{4}<\frac{1}{q}<\frac{d}{2(d+1)}.
	\end{equation*}
	Thus, $(1/r,1/q)\in(N,A)$. At $N$, it follows directly from \eqref{eq: annulus HS}. This completes the proof of \eqref{eq: annulus NA}.
	
	\begin{rem}
		For every $d\geq2$, a restricted type annulus estimate at $A$ \begin{equation*}
			\norm{\sum_j\lambda_j\left|e^{it|\nabla|^{\mu}}|\nabla|^{-s}Qf_j\right|^2}_{L_t^{(d+1)/d}L_x^{(d+1)/(d-1)}}\lesssim\norm{\lambda}_{\ell^{(d+1)/d,1}}
		\end{equation*}follows directly from \cite[Theorem 1.4]{HT}.  
	\end{rem}

	\section{Globalisation}\label{sec: global}
	
	In this section, we use the potential duality principle. Let $d\geq3$ and fix $\mu\neq0,1$. Take an arbitrary $L^2$ orthonormal system $(f_j)_{j\in\mathbb{Z}}$, and let $Q$ be defined as before. Choose the Littlewood--Paley cutoff $\chi\in C_c^\infty$ and define $Q_j$ with $\widehat{Q_jf}(\xi)=\chi(2^{-j}\xi)\widehat{f}(\xi)$ such that $\sum_jQ_j^2=1$. The existence of these operators can be easily justified. 
	
	Denote $T_j(V)=B_{V,Q_j}^{s_\mu}$. 
	
	\begin{lem}
		Let $d\geq3$. Let $(1/r,1/q)\in (OCDA)^o$. Then, \begin{equation}\label{eq: potential sum}
			\left(\sum_{j\in\mathbb{Z}}\norm{T_j(V)}_{\mathfrak{S}^{a}}^{a}\right)^{1/a}\lesssim\norm{V}_{L_t^aL_x^b}
		\end{equation}holds for every $V\in L_t^aL_x^b$, where $a=(q/2)'$, $b=(r/2)'$.
	\end{lem}
	
	\begin{proof}
		For $1<a<\infty$, we define \begin{equation*}
			b_m:=\begin{cases}
				da/(d+a-1),&1<a\leq d+1,\\ da/(2a-2),&d+1<a<\infty.
			\end{cases}
		\end{equation*}Note that for every $1<c<b_m$, the scaling argument and Theorem \ref{thm: annulus} give\begin{equation}\label{eq: T_j}
			\norm{T_j(V)}_{\mathfrak{S}^a}\lesssim2^{dj/c-\mu j(1-1/a)-2s_\mu j}\norm{V}_{L_t^aL_x^c}= 2^{dj(1/c-1/b)}\norm{V}_{L_t^aL_x^c}.
		\end{equation}
		
		Fix $a,b$, and choose $1<c_0<b<c_1<b_m$. Put \begin{equation*}
			\delta_0=d\left(\frac{1}{c_0}-\frac{1}{b}\right),\qquad\delta_1=d\left(\frac{1}{c_1}-\frac{1}{b}\right).
		\end{equation*}Then, $\delta_0>0>\delta_1$.
		
		For $\Gamma=(\Gamma_j)_j\in\ell_{\delta_0}^1(\mathfrak S^{a'})\cap\ell_{\delta_1}^1(\mathfrak S^{a'})$, define a bilinear form \begin{equation}\label{eq: LJ}
			L(V,\Gamma):=\sum_{j\in\mathbb{Z}}\Tr\left(T_j(V)\Gamma_j\right).
		\end{equation}Then, the triangle inequality, the Schatten Hölder inequality and \eqref{eq: T_j} imply \begin{align*}
			|L(V,\Gamma)|&\leq\sum_{j\in\mathbb{Z}}\norm{T_j(V)}_{\mathfrak{S}^a}\norm{\Gamma_j}_{\mathfrak{S}^{a'}}\\&\lesssim \norm{V}_{L_t^aL_x^{c_k}}\norm{\Gamma}_{\ell_{\delta_k}^1(\mathfrak{S}^{a'})},
		\end{align*}where $k=0,1$. Thus, the bilinear real interpolation theorem \cite[Section 3.13, Exercise 5a]{Bergh} provides the following boundedness properties:\begin{equation*}
			L:(L_t^aL_x^{c_0},L_t^aL_x^{c_1})_{\vartheta,p_1}\times(\ell_{\delta_0}^1(\mathfrak{S}^{a'}),\ell_{\delta_1}^1(\mathfrak{S}^{a'}))_{\vartheta,p_2}\longrightarrow(\mathbb{C},\mathbb{C})_{\vartheta,p},
		\end{equation*}where \begin{equation*}
			0<\vartheta<1,\qquad1+\frac{1}{p}=\frac{1}{p_1}+\frac{1}{p_2},\qquad 1\leq p_1,p_2,p\leq\infty.
		\end{equation*}Taking $p_1=a$, $p_2=a'$, choosing $0=(1-\vartheta)\delta_0+\vartheta\delta_1$, using \cite[Theorem 5.6.1]{Bergh} and Lions--Peetre's identity then yield \begin{align*}
			|L(V,\Gamma)|&\lesssim\norm{V}_{(L_t^aL_x^{c_0},L_t^aL_x^{c_1})_{\vartheta,a}}\norm{\Gamma}_{(\ell_{\delta_0}^1,\ell_{\delta_1}^1)_{\vartheta,a'}(\mathfrak{S}^{a'})}\\&\lesssim\norm{V}_{L_t^aL_x^{b,a}}\norm{\Gamma}_{\ell_0^{a'}(\mathfrak{S}^{a'})}.
		\end{align*}Since $b<b_m<a$, we have $L^b=L^{b,b}\hookrightarrow L^{b,a}$. Therefore, \begin{equation*}
			\left|\sum_{j\in\mathbb{Z}}\Tr\left(T_j(V)\Gamma_j\right)\right|\lesssim\norm{V}_{L_t^aL_x^{b}}\norm{\Gamma}_{\ell^{a'}(\mathfrak{S}^{a'})}.
		\end{equation*} Schatten duality \cite[Theorem 3.2]{Simon} and scalar sequence duality \cite[Chapter 1, Lemma 4.2(ii)]{M.Stein 1} prove \eqref{eq: potential sum}.
	\end{proof}
	
	Regarding the bilinear form $L$ in \eqref{eq: LJ}, if we use the sequence Hölder inequality instead of the triangle inequality, then \eqref{eq: potential sum} implies \begin{align*}
		|L(V,\Gamma)|&\leq\left(\sum_{j\in\mathbb{Z}}\norm{B_{V,Q_j}^{s_\mu}}_{\mathfrak{S}^{(q/2)'}}^{(q/2)'}\right)^{1-2/q}\left(\sum_{j\in\mathbb{Z}}\norm{\Gamma_j}_{\mathfrak{S}^{q/2}}^{q/2}\right)^{2/q}\\&\lesssim\norm{V}_{L_t^{(q/2)'}L_x^{(r/2)'}}\left(\sum_{j\in\mathbb{Z}}\norm{\Gamma_j}_{\mathfrak{S}^{q/2}}^{q/2}\right)^{2/q}
	\end{align*}for every $(1/r,1/q)\in(OCDA)^o$. Set $\Gamma_j=\gamma_j:=P_j\gamma P_j$, where $P_j$ is an enlarged cutoff of $Q_j$, meaning that $0\leq P_j\leq 1$, $P_jQ_j=Q_j$, and $\sum_j P_j^2\lesssim 1$. Using \eqref{eq: Vrho} and the triangle inequality yields\begin{align}
		\left|\int_{\mathbb{R}^{d+1}}V(t,x)\sum_{j\in\mathbb{Z}}\rho_{\gamma_{Q_j}(t)}(x)dxdt\right|&\leq \sum_{j\in\mathbb{Z}}\left|\Tr\left(B_{V,Q_j}^{s_\mu}\gamma_{j}\right)\right|\notag\\&\lesssim\norm{V}_{L_t^{(q/2)'}L_x^{(r/2)'}}\left(\sum_{j\in\mathbb{Z}}\norm{\gamma_{j}}_{\mathfrak{S}^{q/2}}^{q/2}\right)^{2/q}.\label{eq: pre}
	\end{align}To approach \eqref{eq: Bv}, we use Littlewood--Paley theory to sum up the pieces of the density operators.
	
	Firstly, we prove the summability of the Schatten norm sequence for positive operators.
	
	\begin{lem}
		For $1\leq\alpha<\infty$ and a positive operator $\Gamma\in\mathfrak{S}^\alpha$, the following statements hold:\begin{itemize}
			\item[1.] For $\norm{u}_{L^2}\leq1$, \begin{equation}\label{eq: scalar Jensen}
				\inpdt{\Gamma u,u}^\alpha\leq\inpdt{\Gamma^\alpha u,u}.
			\end{equation}
			\item[2.] For the Littlewood--Paley cutoffs $Q_j$, \begin{equation}\label{eq: pinch}
				\left(\sum_{j\in\mathbb{Z}}\norm{Q_j\Gamma Q_j}_{\mathfrak{S}^\alpha}^\alpha\right)^{1/\alpha}\leq\norm{\Gamma}_{\mathfrak{S}^\alpha}.
			\end{equation}
		\end{itemize}
	\end{lem} 
	
	\begin{proof}
		Inequality \eqref{eq: scalar Jensen} is trivial if $\norm{u}_{L^2}=0$ or $\alpha=1$. When $0\neq\norm{u}_{L^2}\leq1$ with $\alpha>1$, \eqref{eq: scalar Jensen} follows immediately from \cite[Page 184, (ii)]{Fujii}. 
		
		Fix $j$, and let $(\beta_{j,k},e_{j,k})$ be the nonzero
		eigenvalues and orthonormal eigenvectors of the positive compact
		operator $Q_j\Gamma Q_j$. Then,
		\begin{equation*}
			\beta_{j,k}=\inpdt{\Gamma Q_je_{j,k},Q_je_{j,k}}.
		\end{equation*}
		Since $\norm{Q_je_{j,k}}_{L^2}\leq1$, applying \eqref{eq: scalar Jensen} gives
		\begin{equation*}
			\beta_{j,k}^\alpha
			\leq\inpdt{\Gamma^\alpha Q_je_{j,k},Q_je_{j,k}}.
		\end{equation*}
		Using positivity and trace cyclicity then gives
		\begin{align*}
			\norm{Q_j\Gamma Q_j}_{\mathfrak S^\alpha}^{\alpha}
			&=\sum_k\beta_{j,k}^\alpha\\
			&\leq\sum_k\langle Q_j\Gamma^\alpha Q_je_{j,k},e_{j,k}\rangle\\
			&\leq\Tr(Q_j\Gamma^\alpha Q_j)\\
			&=\Tr
			(\Gamma^{\alpha/2}Q_j^2\Gamma^{\alpha/2}).
		\end{align*}Summing over $j\in\mathbb{Z}$ and using positivity, we obtain
		\begin{align*}
			\sum_{j\in\mathbb{Z}}
			\norm{Q_j\Gamma Q_j}_{\mathfrak S^\alpha}^{\alpha}
			&\le\Tr\left(
			\Gamma^{\alpha/2}
			\left(\sum_{j\in\mathbb{Z}}Q_j^2\right)
			\Gamma^{\alpha/2}\right)\\
			&\le\Tr(\Gamma^\alpha)\\
			&=\norm{\Gamma}_{\mathfrak S^\alpha}^{\alpha}.
		\end{align*}
		This proves \eqref{eq: pinch}.
	\end{proof}
	
	The same proof gives $\sum_j\norm{P_j\Gamma P_j}_{\mathfrak S^\alpha}^\alpha\lesssim\norm{\Gamma}_{\mathfrak S^\alpha}^\alpha$, since $\sum_jP_j^2\lesssim1$. For positive $\gamma$, this bound, standard duality and \eqref{eq: pre} yield \begin{equation}\label{eq: matrix sum}
		\norm{\sum_{j\in\mathbb{Z}}\rho_{\gamma_{Q_j}(t)}}_{L_t^{q/2}L_x^{r/2}}\lesssim\norm{\gamma}_{\mathfrak{S}^{q/2}}.
	\end{equation}It remains to prove, for positive $\gamma$, that if we further have \begin{equation}\label{eq: Littlewood Paley}
		\norm{\rho_{\gamma(t)}}_{L_x^{r/2}}\lesssim\norm{\sum_{j\in \mathbb{Z}}\rho_{\gamma_{Q_j}(t)}}_{L_x^{r/2}},
	\end{equation}where $\gamma(t):=e^{it|\nabla|^\mu}|\nabla|^{-s_\mu}\gamma|\nabla|^{-s_\mu} e^{-it|\nabla|^\mu}$, then using \eqref{eq: matrix sum} we conclude\begin{align*}
		\norm{\sum_{j\in \mathbb{Z}} \lambda_{j}\left|e^{it|\nabla|^\mu}|\nabla|^{-s_\mu}f_j\right|^2}_{L_t^{q/2}L_x^{r/2}}&\leq\norm{\sum_{j\in \mathbb{Z}} |\lambda_{j}|\left|e^{it|\nabla|^\mu}|\nabla|^{-s_\mu}f_j\right|^2}_{L_t^{q/2}L_x^{r/2}}\\&=\norm{\rho_{\dot{\gamma}(t)}}_{L_t^{q/2}L_x^{r/2}}\\&\lesssim\norm{\lambda}_{\ell^{q/2}}
	\end{align*}for every $(\lambda_j)_j$ and every $(1/r,1/q)\in (OCDA)^o$. Here, $\dot{\gamma}(t):=e^{it|\nabla|^\mu}|\nabla|^{-s_\mu}\dot{\gamma}|\nabla|^{-s_\mu} e^{-it|\nabla|^\mu}$, and \begin{equation*}
		\dot{\gamma}:=\sum_{j}|\lambda_j| |f_j\rangle\langle f_j|
	\end{equation*}is positive. By density, this proves Theorem \ref{thm: main} as desired. 
	
	\begin{proof}[Proof of \eqref{eq: Littlewood Paley}] 
		Let $\lambda_j\geq0$. Define the vector-valued function\begin{equation*}
			G_t(x):=\left(\sqrt{\lambda_j}e^{it|\nabla|^\mu}|\nabla|^{-s_\mu}f_j(x)\right)_j.
		\end{equation*}Then, \begin{equation}\label{eq: G}
			\rho_{\gamma(t)}(x)=\norm{G_t(x)}_{\ell^2}^2,\qquad\rho_{\gamma_{Q_j}(t)}(x)=\norm{Q_jG_t(x)}_{\ell^2}^2.
		\end{equation}Since $Q$ is Fourier supported away from the origin, and $r'>1$, the Littlewood--Paley inequality \cite[Proposition 6.1.4]{Grafakos} implies \begin{equation*}
			\norm{\left(\sum_{j\in \mathbb{Z}}\norm{Q_jH}_{\ell^2}^2\right)^{1/2}}_{L_x^{r'}}\lesssim_{d,r,\chi}\norm{H}_{L_x^{r'}(\ell^2)}
		\end{equation*}for every vector-valued function $H\in L_x^{r'}(\ell^2)$. Thus, using Hölder's inequality yields\begin{align*}
			\left|\int_{\mathbb{R}^d}\inpdt{G_t(x),H(x)}_{\ell^2}dx\right|&=\left|\sum_{j\in \mathbb{Z}}\int_{\mathbb{R}^d}\inpdt{Q_jG_t(x),Q_jH(x)}_{\ell^2}dx\right|\\&\leq\norm{\left(\sum_{j\in \mathbb{Z}}\norm{Q_jG_t}_{\ell^2}^2\right)^{1/2}}_{L_x^r}\norm{\left(\sum_{j\in \mathbb{Z}}\norm{Q_jH}_{\ell^2}^2\right)^{1/2}}_{L_x^{r'}}\\&\lesssim\norm{\sum_{j\in \mathbb{Z}}\norm{Q_jG_t}_{\ell^2}^2}_{L_x^{r/2}}^{1/2}\norm{H}_{L_x^{r'}(\ell^2)}.
		\end{align*}Taking the supremum over $\norm{H}_{L_x^{r'}(\ell^2)}\leq1$ and using standard duality gives\begin{equation}\label{eq: LP result}
			\norm{G_t}_{L_x^{r}(\ell^2)}^2\lesssim\norm{\sum_{j\in \mathbb{Z}}\norm{Q_jG_t}_{\ell^2}^2}_{L_x^{r/2}}.
		\end{equation}Substituting \eqref{eq: G} then proves \eqref{eq: Littlewood Paley}.
	\end{proof}
	
	\section{The spatial critical case}
	
	To test \eqref{eq: ortho} on the segment $(O,C)$, we make the initial functions orthogonal, while making their evolved densities contribute at the same spatial point throughout the same time interval. A counterexample on $(O,C)$ then comes from the observation that a large density does not occur only at one isolated time. This gives strong type and weak time type failures immediately. 
	
	\noindent\textbf{The special case. }To see the essential part of the counterexample, we first establish it at $(q,r)=(4,\infty)$. Let us consider the special case $d=2$, $\mu=2$, so that $s=1/2$. Write
	\[
	u_j(t,x)=e^{-it\Delta}|\nabla|^{-1/2}f_j(x).
	\]
	
	Let
	\[
	E_j=\{\xi\in\mathbb R^2:j<|\xi|^2<j+1\},\qquad j\geq1,
	\]
	and define
	\[
	\widehat f_j(\xi)=2\sqrt{\pi}\,1_{E_j}(\xi).
	\]
	The set \(E_j\) is the annulus between the circles of radii \(\sqrt j\) and \(\sqrt{j+1}\). Its area is
	\[
	|E_j|
	=\pi(\sqrt{j+1})^2-\pi(\sqrt j)^2
	=\pi.
	\]
	By Plancherel,
	\[
	\begin{aligned}
		\|f_j\|_{L^2}^2
		&=(2\pi)^{-2}\int_{\mathbb R^2}|\widehat f_j(\xi)|^2\,d\xi\\
		&=\frac1{4\pi^2}\int_{E_j}4\pi\,d\xi\\
		&=\frac1{4\pi^2}\cdot4\pi\cdot\pi\\
		&=1.
	\end{aligned}
	\]
	For \(j\ne k\), the annuli are disjoint, so
	\[
	\begin{aligned}
		\langle f_j,f_k\rangle
		=(2\pi)^{-2}
		\int_{\mathbb R^2}\widehat f_j(\xi)
		\overline{\widehat f_k(\xi)}\,d\xi=0.
	\end{aligned}
	\]
	Therefore \((f_j)\) is orthonormal.
	
	Fix $j\in\mathbb{N}$ and write $r=|\xi|$, $\theta=\arg\xi$, $u=r^2$, and $v=u-j$ to obtain
	\begin{align*}
		u_j(t,0)
		&=(2\pi)^{-2}
		\int_{\mathbb R^2}
		e^{it|\xi|^2}|\xi|^{-1/2}\widehat f_j(\xi)\,d\xi\\&=\frac{2\sqrt\pi}{(2\pi)^2}
		\int_{E_j}e^{it|\xi|^2}|\xi|^{-1/2}\,d\xi\\&=\frac{2\sqrt\pi}{4\pi^2}
		\int_0^{2\pi}\int_{\sqrt j}^{\sqrt{j+1}}
		e^{itr^2}r^{-1/2}r\,dr\,d\theta\\
		&=\frac1{\sqrt\pi}
		\int_{\sqrt j}^{\sqrt{j+1}}e^{itr^2}r^{1/2}\,dr\\&=\frac1{2\sqrt\pi}
		\int_j^{j+1}e^{itu}u^{-1/4}\,du\\&=\frac{e^{ijt}}{2\sqrt\pi}
		\int_0^1e^{itv}(j+v)^{-1/4}\,dv.
	\end{align*}
	Thus,
	\[
	|u_j(t,0)|
	=\frac1{2\sqrt\pi}
	\left|\int_0^1e^{itv}(j+v)^{-1/4}\,dv\right|.
	\]
	
	We derive a lower bound for the absolute
	value of $u_j(t,0)$ on the interval $|t|\leq1$. In this case, \(0\le v\le1\),
	\[
	|tv|\le1,
	\qquad
	\cos(tv)\ge\cos(1)>0.
	\]
	Using the inequality \(|z|\ge\operatorname{Re}z\),
	\begin{align}
		|u_j(t,0)|
		&\ge\frac1{2\sqrt\pi}
		\operatorname{Re}\int_0^1e^{itv}(j+v)^{-1/4}\,dv\notag\\
		&=\frac1{2\sqrt\pi}
		\int_0^1\cos(tv)(j+v)^{-1/4}\,dv\notag\\
		&\ge\frac{\cos(1)}{2\sqrt\pi}
		\int_0^1(j+v)^{-1/4}\,dv.\label{eq: uj}
	\end{align}
	Since \(j\ge1\),
	\[
	j+v\le j+1\le2j.
	\]
	So
	\[
	(j+v)^{-1/4}\ge(2j)^{-1/4}.
	\]
	Consequently, we obtain
	\begin{equation*}
		|u_j(t,0)|^2\ge c_0j^{-1/2},
		\qquad |t|\le1\textcolor{red}{,}
	\end{equation*}
	where the fixed constant is
	\[
	c_0=\frac{\cos^2(1)}{4\sqrt2\pi}>0.
	\]
	
	Take
	\[
	\lambda_j=j^{-1/2},\qquad 1\le j\le N,
	\]
	and write
	\[
	H_N=\sum_{j=1}^N\frac1j.
	\]
	Since every coefficient is positive, \eqref{eq: uj} gives
	\[
	\begin{aligned}
		\sum_{j=1}^N\lambda_j|u_j(t,0)|^2
		\ge c_0\sum_{j=1}^N j^{-1/2}j^{-1/2}=c_0H_N,
		\qquad |t|\le1,
	\end{aligned}
	\]with $\|\lambda\|_{\ell^2}=H_N^{1/2}$. Thus, the density has size at least \(H_N\), while the coefficient norm has size only \(H_N^{1/2}\). We justify this result on a ball near the spatial origin.
	
	For each fixed \(j,t\), 
	\[
	e^{it|\xi|^2}|\xi|^{-1/2}\widehat f_j(\xi)
	\]
	belongs to \(L_\xi^1\). Indeed, it is bounded and supported on the bounded ring \(E_j\).
	Its inverse Fourier transform \(u_j(t,x)\) is therefore continuous in \(x\). Hence the finite sum
	\[
	\rho_N(t,x)=\sum_{j=1}^N\lambda_j|u_j(t,x)|^2
	\]
	is continuous in \(x\).
	For any \(\varepsilon>0\), continuity gives a ball around \(0\) on which
	\[
	\rho_N(t,x)>\rho_N(t,0)-\varepsilon.
	\]
	Therefore
	\[
	\|\rho_N(t,\cdot)\|_{L_x^\infty}
	\ge\rho_N(t,0)-\varepsilon.
	\]
	Letting \(\varepsilon\downarrow0\) gives $\|\rho_N(t,\cdot)\|_{L_x^\infty}
	\ge\rho_N(t,0)$. Thus, for this fixed $|t|\le1$,
	\[
	\|\rho_N(t,\cdot)\|_{L_x^\infty}\ge c_0H_N,
	\]
	and
	\[
	\begin{aligned}
		\|\rho_N\|_{L_t^2L_x^\infty}
		&\ge
		\left(\int_{-1}^1
		\|\rho_N(t,\cdot)\|_{L_x^\infty}^2\,dt\right)^{1/2}\\
		&\ge
		\left(\int_{-1}^1c_0^2H_N^2\,dt\right)^{1/2}\\
		&=\sqrt2\,c_0H_N.
	\end{aligned}
	\]
	Dividing by the coefficient norm gives \begin{equation*}
		\frac{\|\rho_N\|_{L_t^2L_x^\infty}}{\|\lambda\|_{\ell^2}}
		\ge\sqrt2\,c_0H_N^{1/2}.
	\end{equation*}
	Letting $N\to\infty$ proves the strong type failure. 
	
	To see the failure of the weak-time estimate, we put $g_N(t)=\|\rho_N(t,\cdot)\|_{L_x^\infty}$, and then
	\[
	g_N(t)\ge c_0H_N\qquad\text{on }[-1,1].
	\]
	Choose
	\[
	\tau=\frac{c_0H_N}{2}.
	\]
	Then,
	\[
	\|g_N\|_{L^{2,\infty}}
	=\sup_{\tau>0}
	\tau\,|\{t:g_N(t)>\tau\}|^{1/2}\ge\frac{c_0H_N}{\sqrt{2}}.
	\]
	This lower bound gives a logarithmically divergent ratio, proving failure of the weak-time estimate.
	
	To obtain a counterexample with smooth data, choose a fixed function
	\[
	h\in C_c^\infty((0,1)),
	\qquad h\ge0,
	\qquad \int_0^1h(v)^2\,dv=1,
	\]
	and define
	\[
	\widehat f_j(\xi)=2\sqrt\pi\,h(|\xi|^2-j).
	\]Then, $(f_j)$ is orthonormal, and similarly we obtain \begin{equation*}
		|u_j(t,0)|
		\ge
		\frac{\cos(1)}{2\sqrt\pi}(2j)^{-1/4}
		\int_0^1h(v)\,dv.
	\end{equation*}The same argument therefore proves failure of both the strong estimate and the weak-time estimate for these smooth data.
	
	\vspace{2ex}
	
	\noindent\textbf{The general case.}
	Take
	\[
	d\ge2,\qquad \mu\in\mathbb R\setminus\{0,1\},
	\qquad 2<q<\infty,
	\qquad s=\frac d2-\frac\mu q.
	\]
	Let
	\[
	\omega_{d-1}=|\mathbb S^{d-1}|.
	\]
	Choose the same nonnegative smooth \(h\) as above, with \(\|h\|_{L^2}=1\), and define
	\begin{equation*}
		\widehat f_j(\xi)
		=c_{d,\mu}|\xi|^{(\mu-d)/2}
		h(|\xi|^\mu-j),
		\qquad
		c_{d,\mu}=(2\pi)^{d/2}
		\sqrt{\frac{|\mu|}{\omega_{d-1}}},
	\end{equation*}
	where \(|\xi|^{(\mu-d)/2}\) is chosen to make the normalisation simple.
	Indeed, polar coordinates give
	\[
	\begin{aligned}
		\|f_j\|_{L^2}^2
		&=(2\pi)^{-d}c_{d,\mu}^2\omega_{d-1}
		\int_0^\infty
		r^{\mu-d}h(r^\mu-j)^2r^{d-1}\,dr\\
		&=(2\pi)^{-d}c_{d,\mu}^2\omega_{d-1}
		\int_0^\infty h(r^\mu-j)^2r^{\mu-1}\,dr.
	\end{aligned}
	\]
	With \(u=r^\mu\), the change of variables gives
	\[
	\int_0^\infty h(r^\mu-j)^2r^{\mu-1}\,dr
	=\frac1{|\mu|}\int_0^\infty h(u-j)^2\,du
	=\frac1{|\mu|}.
	\]
	Thus, $\|f_j\|_{L^2}=1$. Now compute the solution at the origin for $U_s(t)=e^{it|\nabla|^\mu}|\nabla|^{-s}$:
	\[
	\begin{aligned}
		U_s(t)f_j(0)
		&=(2\pi)^{-d}c_{d,\mu}\omega_{d-1}
		\int_0^\infty
		e^{itr^\mu}
		r^{-s}r^{(\mu-d)/2}
		h(r^\mu-j)r^{d-1}\,dr\\
		&=(2\pi)^{-d}c_{d,\mu}\omega_{d-1}
		\int_0^\infty
		e^{itr^\mu}
		r^{d/2+\mu/2-s-1}h(r^\mu-j)\,dr\\
		&=(2\pi)^{-d/2}\sqrt{\frac{\omega_{d-1}}{|\mu|}}e^{ijt}\int_0^1e^{itv}(j+v)^{1/q-1/2}h(v)\,dv,
	\end{aligned}
	\]
	where \(u=r^\mu\). Since \(q>2\),
	\[
	\eta:=\frac1q-\frac12<0.
	\]
	For \(j\ge1\) and \(0<v<1\),
	\[
	(j+v)^\eta\ge(2j)^\eta.
	\]
	Taking real parts as before gives
	\[
	|U_s(t)f_j(0)|
	\ge (2\pi)^{-d/2}\sqrt{\frac{\omega_{d-1}}{|\mu|}}\cos(1)\,2^\eta j^\eta
	\int_0^1h(v)\,dv,
	\qquad |t|\le1.
	\]
	Therefore, for fixed $|t|\leq1$,\begin{equation*}
		|U_s(t)f_j(0)|^2\ge c\,j^{2/q-1},
	\end{equation*}
	with \(c>0\) independent of \(j\).
	Finally, choose
	\[
	\lambda_j=j^{-2/q},
	\qquad 1\le j\le N.
	\]
	Then
	\[
	\sum_{j=1}^N\lambda_j|U_s(t)f_j(0)|^2
	\ge c\sum_{j=1}^N j^{-2/q}j^{2/q-1}
	=cH_N,
	\]
	whereas
	\[
	\begin{aligned}
		\|\lambda\|_{\ell^{q/2}}
		&=\left(\sum_{j=1}^N
		(j^{-2/q})^{q/2}\right)^{2/q}\\
		&=H_N^{2/q}.
	\end{aligned}
	\]
	Therefore, we obtain
	\begin{equation*}
		\frac{
			\left\|\sum_{j=1}^N\lambda_j|U_s(t)f_j|^2
			\right\|_{L_t^{q/2}L_x^\infty}
		}{
			\|\lambda\|_{\ell^{q/2}}
		}
		\gtrsim H_N^{\,1-2/q}\longrightarrow\infty,
	\end{equation*}
	since \(q>2\). The same lower bound on $[-1,1]$ gives this divergence with $L_t^{q/2,\infty}L_x^\infty$.
	
	\section{Two-dimensional estimates}
	
	Let $d=2$. To prove the annulus estimate on $OCA\setminus[C,A]$, it suffices to show that the inequality\begin{equation}\label{eq: OC ortho annulus}
		\norm{\sum_{j}\lambda_j\left|e^{it|\nabla|^\mu}|\nabla|^{-s_\mu}Qf_j\right|^2}_{L_t^{q/2}L_x^{\infty}(\mathbb{R}\times\mathbb{R}^2)}\lesssim\norm{\lambda}_{\ell^{q/2}}
	\end{equation}holds for every $2<q<\infty$. Here, $s_\mu=1-\mu/q$. Since the argument in Section \ref{sec: annulus} for $[N,A)$ remains valid for $d=2$, we obtain \eqref{eq: annulus NA} on $[O,A)$ for $d=2$. Using a similar complex interpolation argument yields \eqref{eq: annulus NA} on $OCA\setminus[C,A]$, and the summation method in Section \ref{sec: global} applies for $d=2$ with $b_m=\min\{2a/(a+1),a/(a-1)\}$. For $a=(q/2)'$ and $b=(r/2)'$, the region below gives $1<b<b_m<a$. This gives the following statement. \begin{thm}\label{thm: d2}
		Let $\mu\in\mathbb R\setminus\{0,1\}$. Suppose that the admissible pair $(q,r)$ satisfies \begin{equation}\label{geo: region 2d}
			2< q, r<\infty,\qquad q<\frac{r}{2},\qquad\frac{1}{q}+\frac{1}{r}<\frac{1}{2}.
		\end{equation}Then the inequality\begin{equation}\label{eq: main ortho 2d}
			\norm{\sum_j \lambda_j\left|e^{it|\nabla|^\mu}|\nabla|^{-s_\mu}f_j\right|^2}_{L_t^{q/2}L_x^{r/2}(\mathbb{R}\times\mathbb{R}^2)}\lesssim\norm{\lambda}_{\ell^{q/2}}
		\end{equation}holds for every orthonormal system $(f_j)_j\subset L^2(\mathbb{R}^2)$ and every coefficient sequence $\lambda\in\ell^{q/2}$, with \begin{equation}\label{eq: main scaling 2d}
			s_\mu=1-\frac{2}{r}-\frac{\mu}{q}.
		\end{equation}
	\end{thm}
	
	We remark that since $(f_j)_j$ is an orthonormal system, the normalised Fourier transforms $((2\pi)^{-1}\widehat{f_j})_j$ also form an orthonormal system in $L^2$, and thus Bessel's inequality gives  \begin{align*}
		\sum_j\left|e^{it|\nabla|^\mu}|\nabla|^{-s_\mu}Qf_j(x)\right|^2&=\frac{1}{16\pi^4}\sum_j\left|\inpdt{\widehat{f_j}(\xi),e^{-i(x\cdot\xi+t|\xi|^\mu)}|\xi|^{-s_\mu}\chi(\xi)}_{L_\xi^2}\right|^2\\&\leq\frac{1}{4\pi^2}\norm{e^{-i(x\cdot\xi+t|\xi|^\mu)}|\xi|^{-s_\mu}\chi(\xi)}_{L_\xi^2}^2\\&\lesssim\norm{\chi}_{L_\xi^2}^2.
	\end{align*}Then, using Hölder's inequality gives\begin{align*}
		\sum_j|\lambda_j|\left|e^{it|\nabla|^\mu}|\nabla|^{-s_\mu}Qf_j(x)\right|^2&=\sum_j|\lambda_j|\left|e^{it|\nabla|^\mu}|\nabla|^{-s_\mu}Qf_j(x)\right|^{2/\beta}\left|e^{it|\nabla|^\mu}|\nabla|^{-s_\mu}Qf_j(x)\right|^{2/\beta'}\\&\lesssim_{\beta,\chi}\left(\sum_j|\lambda_j|^{\beta}\left|e^{it|\nabla|^\mu}|\nabla|^{-s_\mu}Qf_j(x)\right|^{2}\right)^{1/\beta},
	\end{align*}where $1\leq\beta\leq q/2$. The triangle inequality thus implies \begin{equation}\label{eq: bessel result}
		\norm{\sum_j\lambda_j\left|e^{it|\nabla|^\mu}|\nabla|^{-s_\mu}Qf_j\right|^2}_{L_t^{q/2}L_x^\infty}\lesssim_{\beta,\chi}\left(\sum_j|\lambda_j|^{\beta}\norm{e^{it|\nabla|^\mu}|\nabla|^{-s_\mu}Qf_j}_{L_t^{q/\beta}L_x^\infty}^2\right)^{1/\beta}.
	\end{equation}The right hand side of \eqref{eq: bessel result} requires the scalar annulus Strichartz estimate\begin{equation}\label{eq: OC scalar}
		\norm{e^{it|\nabla|^\mu}|\nabla|^{-s_\mu}Qf}_{L_t^{p}L_x^\infty}\lesssim_{\chi}\norm{f}_{L^2}
	\end{equation}for some $2\leq p\leq\infty$. At $p=2$ it fails when $\mu=2$; see \cite{Montgomery-Smith} and \cite[Theorem 1.1]{T.Tao 1}. This illustrates the difficulty in deducing \eqref{eq: OC ortho annulus} directly from \eqref{eq: bessel result}, since $q/\beta=2$ if $\beta=q/2$; see also the concluding discussion in \cite[Section 3]{Bez}.
	
	\begin{proof}[Proof of Theorem \ref{thm: d2}]
		To prove \eqref{eq: OC ortho annulus}, we use potential duality. It suffices to show that\begin{equation}\label{eq: OC dual annulus}
			\norm{\int_{\mathbb{R}} U(t)^\ast V(t,x)U(t)dt}_{\mathfrak{S}^{q/(q-2)}}\lesssim_{\mu,\chi,q}\norm{V}_{L_t^{q/(q-2)}L_x^1}
		\end{equation}holds for every $0\leq V\in L_t^{q/(q-2)}L_x^1$, where $U(t)=e^{it|\nabla|^\mu}|\nabla|^{-s_\mu}Q$. We prove \eqref{eq: OC dual annulus} for every integer $p=q/(q-2)\geq2$. With the exponent of $|\nabla|$ in $U(t)$ fixed, interpolation with $O$ then gives every $q>2$, since $2p/(p-1)\downarrow2$.

		Fix an integer $p=q/(q-2)\geq2$. As in \cite[Section 4.1]{Frank1}, without any loss of generality, we may assume that $V\geq0$ is bounded and compactly supported. By cyclicity of the trace, we have  \begin{align*}
			\norm{\int_{\mathbb{R}} U(t)^\ast V(t,x)U(t)dt}_{\mathfrak{S}^p}^p&=\Tr\left(\left(\int_{\mathbb{R}}U(t)^\ast V(t,x)U(t)dt\right)^p\right)\\&=\int_{\mathbb{R}^{3p}}\prod_{j=1}^{p}\left(V(t_j,x_j)K(t_j-t_{j+1},x_j-x_{j+1})\right)dt_1\cdots dt_pdx_1\cdots dx_p, 
		\end{align*}where $t_{p+1}:=t_1$, $x_{p+1}:=x_1$, and\begin{equation*}
			K(t,x):=\frac{1}{4\pi^2}\int_{\mathbb{R}^2}e^{i(x\cdot\xi+t|\xi|^\mu)}|\xi|^{-2s_\mu}\chi(\xi)^2d\xi.
		\end{equation*}Here, $K(t-t',x-y)$ is the kernel of $U(t)U(t')^\ast$. By applying the argument for \eqref{eq: dispersive} with the smooth annulus amplitude $|\xi|^{-2s_\mu}\chi(\xi)^2$, we have \begin{equation}\label{eq: dispersive 2d}
			\sup_x|K(t,x)|\lesssim(1+|t|)^{-1}.
		\end{equation}Thus, \begin{align}
			\norm{\int_{\mathbb{R}} U(t)^\ast V(t,x)U(t)dt}_{\mathfrak{S}^p}^p&\leq	\int_{\mathbb{R}^{3p}}\prod_{j=1}^{p}\left|\left(V(t_j,x_j)K(t_j-t_{j+1},x_j-x_{j+1})\right)\right| dt_1\cdots dt_pdx_1\cdots dx_p\notag\\&\lesssim\int_{\mathbb{R}^{3p}}\prod_{j=1}^{p}\left(V(t_j,x_j)\left(1+|t_j-t_{j+1}|\right)^{-1}\right) dt_1\cdots dt_pdx_1\cdots dx_p\notag\\&=\int_{\mathbb{R}^{p}}\prod_{j=1}^{p}\left(v(t_j)k(t_j-t_{j+1})\right) dt_1\cdots dt_p\notag\\&=:I,\label{eq: I}
		\end{align}where\begin{equation*}
			v(t)=\int_{\mathbb{R}^2}V(t,x)dx,\qquad k(t)=(1+|t|)^{-1}.
		\end{equation*}
		
		Make the change of variables\begin{equation*}
			t:=t_1,\qquad u_j:=t_j-t_1.
		\end{equation*}Then, $u_1=u_{p+1}=0$, and \begin{equation*}
			t_j=t+u_j,\qquad t_j-t_{j+1}=u_j-u_{j+1}.
		\end{equation*}The Jacobian determinant of this change of variables is $1$, so Hölder's inequality and the translation invariance of the Lebesgue norm imply\begin{align}
			I&=\int_{\mathbb{R}^{p-1}}\left(\int_{\mathbb{R}}\prod_{j=1}^pv(t+u_j)dt\right)\prod_{j=1}^pk(u_j-u_{j+1})du_2\cdots du_p\notag\\&\leq\int_{\mathbb{R}^{p-1}}\left(\prod_{j=1}^p\norm{v(\cdot+u_j)}_{L^p}\right)\prod_{j=1}^pk(u_j-u_{j+1})du_2\cdots du_p\notag\\&=\norm{v}_{L^p}^p M,\label{eq: M}
		\end{align}where \begin{equation*}
			M:=\int_{\mathbb{R}^{p-1}}\prod_{j=1}^pk(u_j-u_{j+1})du_2\cdots du_p.
		\end{equation*}
		
		Observe that $M$ has a repeated convolution structure with $p$ factors. Since $u_1=0$, we obtain\begin{equation*}
			M=\int_{\mathbb{R}}k_1(-u_2)(k_2\ast\cdots\ast k_p)(u_2)du_2=(k_1\ast\cdots\ast k_p)(0),
		\end{equation*}where $k_1=\cdots=k_p=k$. Using Young's convolution inequality then yields\begin{equation}\label{eq: M p}
			M\leq\norm{k_1\ast\cdots\ast k_p}_{L^\infty}\leq\prod_{j=1}^p\norm{k_j}_{L^{p'}}=\norm{k}_{L^{q/2}}^p.
		\end{equation}Since $q>2$, \begin{equation*}
			\norm{k}_{L^{q/2}}^{q/2}=\int_{\mathbb{R}}\frac{1}{(1+|t|)^{q/2}}dt<\infty.
		\end{equation*}Therefore, substituting \eqref{eq: M p} into \eqref{eq: M} and \eqref{eq: I} establishes \eqref{eq: OC dual annulus}. The interpolation and globalisation described above then prove \eqref{eq: main ortho 2d}.
	\end{proof}
	
	\vspace{4ex}
	
	\textbf{Contact.} Yinghao Gong, Nagoya University, \textit{Email: gong.yinghao.m7@s.mail.nagoya-u.ac.jp}.

\end{document}